%% file: ORC-Prediction.tex
\documentclass[letterpaper, 10 pt, conference]{ieeeconf}  

\IEEEoverridecommandlockouts     

\usepackage{color,amsmath,amstext,amsfonts,amssymb, mathrsfs,mathtools}
\usepackage{xcolor, graphicx,algorithmic,algorithm}
\usepackage[algo2e,linesnumbered]{algorithm2e}
\usepackage{tikz}
\usepackage{setspace}
\usepackage{array} 
\usepackage{multirow}
\usepackage{booktabs}  
\usepackage{bbm}
\usepackage{wasysym}
\usepackage{textcomp}
\usepackage{hyperref} 
\usepackage[sort,compress]{cite}

\newtheorem{assumption}{Assumption}

\newtheorem{theorem}{Theorem}

\newtheorem{remark}{Remark}

\renewcommand{\thesubsubsection}{{\it \Alph{subsection}.\arabic{subsubsection}.}}

\SetKwInput{KwInput}{Input}

\newcounter{l1}
\newcounter{l2}
\newcounter{l3}
\newcommand{\bdotlist}{\begin{list}{$\bullet$}{}}
\newcommand{\bboxlist}{\begin{list}{$\Box$}{}}
\newcommand{\bbboxlist}{\begin{list}{\raisebox{.005in}{{\tiny $\blacksquare$ \ \ }}}{}}
\newcommand{\bdashlist}{\begin{list}{$-$}{} }
\newcommand{\blist}{\begin{list}{}{} }
\newcommand{\barablist}{\begin{list}{\arabic{l1}}{\usecounter{l1}}}
\newcommand{\balphlist}{\begin{list}{(\alph{l2})}{\usecounter{l2}}}
\newcommand{\bAlphlist}{\begin{list}{\Alph{l2}.}{\usecounter{l2}}}
\newcommand{\bdiamlist}{\begin{list}{$\diamond$}{}}
\newcommand{\bromalist}{\begin{list}{(\roman{l3})}{\usecounter{l3}}}

\providecommand{\norm}[1]{\lVert#1\rVert}

\newcommand{\beq}{\begin{equation}}
\newcommand{\eeq}{\end{equation}}

\allowdisplaybreaks

\title{Online Robust Control of Linear Dynamical Systems with Limited Prediction}

\author{Deepan Muthirayan, Dileep Kalathil, and Pramod P. Khargonekar
\thanks{This work is supported in part by the National Science Foundation under Grant ECCS-1839429.
D. Muthirayan and P. P. Khargonekar are with the Department of Electrical Engineering and Computer Sciences, University of California Irvine, Irvine, CA (emails: deepan.m@uci.edu, pramod.khargonekar@uci.edu). Dileep Kalathil is with the Department of Electrical and Computer Engineering, Texas A\& M University (email:dileep.kalathil@tamu.edu).}
}

\begin{document}

\maketitle
\thispagestyle{empty}
\pagestyle{empty}

\input{abstract}

\input{introduction}

\input{formulation}

\input{relatedworks} 

\input{algorithm-results} 

\input{analysis} 

\input{conclusion} 

\bibliographystyle{IEEEtran} 
\bibliography{Refs.bib}

\input{appendix} 

\end{document}

%% file: abstract.tex
\begin{abstract}
We study the online robust control problem for linear dynamical systems with disturbances and uncertainties in the cost functions, with limited preview of the future disturbances and the cost functions, $N$. 
Our goal is to find an online control policy that can minimize the disturbance gain, defined as the ratio of the cumulative cost and the cumulative energy in the disturbances over a period of time, in the face of the uncertainties, and characterize its achievable gain in terms of the system relevant parameters. Our goals contrast with prior online control works for the same problem, which either focus on minimizing the static regret, a weaker performance metric, or assume a very large preview of the future uncertainties.  
Specifically, we consider a class of cost functions characterized by $\beta$ ($\beta < 1$), a number whose inverse bounds the variation of the cost functions. 
We propose a novel variation of the Receding Horizon Control as the online control policy. We show that, under standard system assumptions, when $N > 4/\beta^3$, the proposed algorithm can achieve a disturbance gain $(2/\beta+\rho(N)) \overline{\gamma}^2$, where $\overline{\gamma}^2$ is the best (minimum) possible disturbance gain for an oracle policy with full knowledge of the cost functions and disturbances, with $\rho(N) = O(1/N)$. We also demonstrate through simulations that the proposed policy satisfies the derived bounds and is consistently better than the standard RHC approach. 

\end{abstract}

%% file: introduction.tex
\section{Introduction}

The control of dynamical systems with uncertainties such as modeling errors, parametric uncertainty, and disturbances is a central challenge in control theory. There is vast literature in the field on control synthesis for systems with such uncertainties. The \textit{robust control} literature studies the problem of feedback control design for stability and performance guarantees with modeling uncertainty and disturbances; see \cite{skogestad2007multivariable}. The \textit{adaptive control} literature studies the control of systems with parametric uncertainty; see \cite{sastry2011adaptive, ioannou2012robust, aastrom2013adaptive}. Typically, these approaches are concerned with stability and asymptotic performance guarantees. In many practical control applications though, the control task is of a finite duration and therefore the performance over a finite period of time becomes crucially important. This is the central focus in the field of online control, which is developing online adaptive algorithms that optimize the performance over a finite period of time. 

A standard formulation in online control considers systems with uncertainties that can be quite general like unknown system parameters, disturbances or even unknown and possibly time varying cost functions. The controller at any time step can observe only the state of the system and the cost of its action (bandit feedback) or the full cost function for that time step (full information feedback). The challenge in online control is learning from the incremental information acquired from every time step and simultaneously optimizing the performance. While the challenge of adaptation is the same as in classical adaptive control, a key difference in online control is the finite performance metric such as the regret. A typical regret metric is defined as the difference between the cumulative cost of the online controller and the cumulative cost of the best policy from a comparator class over a given period of time. The non-asymptotic performance objective is the key challenge compared to the classical adaptive control literature, and typically, solving online control problems requires combining techniques from statistical learning, online optimization and control. 

In this work we study the online control of linear dynamical systems with disturbances. Most works that study this problem consider the regret metric as the objective with linear feedback policies as the comparator class \cite{dean2018regret, cohen2019learning, mania2019certainty, agarwal2019online, agarwal2019logarithmic, simchowitz2020improper}. However, the absolute performance of the online controller designed to minimize the regret metric crucially depends on the specific comparator class that is used. While the class of linear feedback policies is the optimal for the Linear Quadratic Regulator (LQR) problem, this need not hold true for more generic cost functions such as convex functions, such as the cost functions considered in the online control problems in \cite{ agarwal2019online, agarwal2019logarithmic, simchowitz2020improper}. In this work our goal is to find an online control policy that {\it minimizes the disturbance gain, defined as the ratio of the cumulative cost and the cumulative energy in the disturbances over a given period of time.} Unlike the regret metric, which is a relative metric, disturbance gain is a widely used measure in control for characterizing the absolute performance.

Our problem is similar to the well studied $\mathcal{H}_\infty$ robust control problem. Here, the objective is to minimize the worst-case gain of the energy in the disturbance to the energy in the output; see \cite{doyle1988state}. This problem has been well studied for linear dynamic systems and quadratic costs. Explicit characterization of the optimal state feedback controller has been given for this problem; see \cite{doyle1988state}. Such characterization is harder when the cost functions are general and even harder in the online setting. In such cases, \textit{how do we develop online control algorithms to minimize the gain over the disturbances over a finite period of time}?

In the online setting we study, the system is a known linear dynamic system with arbitrary disturbance in its dynamics. The cost functions are arbitrary, time varying and unknown apriori. Such settings are relevant in many applications like, for eg., Autonomous Vehicles (AV), where the cost functions are dependent on conditions like traffic which can be unpredictable. We specifically study the setting where the controller has limited preview of the future cost functions and disturbances. This is reasonable considering that in many practical problems such as robotics \cite{baca2018model, shi2019neural}, energy systems \cite{vazquez2016model}, data-center management \cite{lazic2018data} etc., a fixed finite-length preview of the future cost functions and/or disturbances are available to compute the control input. Specifically we address the following questions: \textit{how do we use the limited preview of future cost functions and disturbances to adapt online to minimize the overall gain over the disturbances}? \textit{What is the achievable gain in terms of the system or environment parameters}? Specifically we address the question of what is the achievable gain in terms of the preview length and the variability in the cost functions.

\subsection{Main Contribution}

We present online control algorithms for linear dynamical systems with arbitrary and time varying cost functions and arbitrary disturbances. Specifically, we address the problem of disturbance gain minimization over a finite time with limited preview of future cost functions and disturbances, $N$. We consider a class of cost functions characterized by $\beta$ ($\beta < 1$), a number whose inverse bounds the variation of the cost functions. The cost functions that satisfy this condition are sufficiently general and include the time varying LQR. Thus our guarantees are applicable to a general class of systems. We present a novel extension of the standard Receding Horizon Controller (RHC) as the online control policy and characterize its disturbance gain in terms of the system relevant parameters, i.e., preview length and $\beta$. To the best of our knowledge ours is the first work to give such a characterization of the disturbance gain for any type of predictive control. 

Specifically, we show that, when the initial state's cost is zero, and $N > 4\zeta^3$ for some $\zeta \geq 1/\beta$, for any $T > N$, the proposed RHC's total cost $J$ is bounded by \beq J \leq (2\zeta+\rho(N))\overline{\gamma}^2 \sum_{t=1}^T \norm{w_t}^2, \nonumber \eeq  with $\rho(N) = \mathcal{O}(1/N)$. Thus, the disturbance gain achieved by the proposed RHC is given by $(2\zeta+\rho(N))\overline{\gamma}^2$. Thus, we show that {\it the online RHC's disturbance gain differs from $\overline{\gamma}^2$ only by a system dependent constant factor and reduces with preview $N$ as $\mathcal{O}(1/N)$}. We also present numerical examples to demonstrate the effectiveness of our proposed algorithm over the standard RHC.


\subsection{Comparison to Relevant Literature}

Receding Horizon Control, also known as Model Predictive Control (MPC), address the class of problems where a certain preview of future cost functions are available to compute the control input. The RHC is a well studied methodology in the control literature \cite{mayne2000constrained, angeli2011average, camacho2013model, grune2014asymptotic, angeli2016theoretical, borrelli2017predictive, rosolia2017learning, baca2018model, grune2020economic}. However, the control theory literature on RHC mostly focus on asymptotic performance guarantees. Works that present finite time guarantees focus on Input to State/Output Stabilization (ISS) \cite{lazar2008input, lale2021model} or set convergence in finite time \cite{anderson2018finite}. Thus, these works address specific cost conditions and do not consider the variability in cost functions like we do. 

In the machine learning literature, there are works that study the problem of online control for general cost functions with preview \cite{li2019online, goel2020regret, yu2020power, lin2021perturbation}. The metric that these works study is a stronger notion of regret called the {\it dynamic regret}, which is the regret where the comparator class is not restricted to any policy class. Thus, dynamic regret is a stronger metric than the static regret and is a measure of the absolute performance. In contrast to \cite{li2019online, goel2020regret, yu2020power, lin2021perturbation}, we characterize the disturbance gain. In addition, \cite{li2019online, goel2020regret, yu2020power, lin2021perturbation} assume large preview or full knowledge of the future cost functions, while we study the setting where {\it the preview of the future cost functions is limited}. 


\subsection{Notations}

For a sequence (of vectors/function) $(a_{1}, a_2,..., a_T)$, we denote $a_{1:t} = (a_{1}, a_2,..., a_T)$. For a matrix $M$, we denote its transpose by $M^\top$. We denote the maximum eigenvalue of a matrix $M$ by $\lambda_{\text{max}}(M)$. The two norm of a vector is denoted by $\norm{.}_2$. When two matrices $M_1$ and $M_2$ are related by $M_1 \geq M_2$ then it implies that $M_1 - M_2$ is positive semi-definite. Similarly, when $M_1 > M_2$ it implies that $M_1 - M_2$ is positive definite. We denote $\mathbb{R}^n$ as the $n-$dimensional euclidean space and $\mathbb{R}_{+}$ as the non-negative part of the real line. 
We compactly denote a time interval from $t_1$ to $t_2$ by $[t_1,t_2]$.

%% file: formulation.tex
\section{Problem Preliminaries}
\label{sec:prob-prelim}
We consider the online control problem of a linear dynamical system with cost functions and disturbance that are arbitrary. In particular, the system evolution is represented by the equation 
\begin{equation}
x_{t+1} = Ax_t+Bu_t+w_t,
\label{eq:stateequation}
\end{equation} 
where $x_{t} \in \mathbb{R}^n$ is the system state, $u_t \in \mathbb{R}^m$ is the control input, and $w_t \in \mathbb{R}^n$ is the \textit{arbitrary} disturbance at time step $t$.  At each time step $t$, the control algorithm selects a control input $u_{t}$ and incurs a cost $c_{t}(x_{t}, u_{t})$. The cost function sequence $\{c_{t}(\cdot, \cdot)\}$ is assumed to be \textit{arbitrary} and hence unknown to the algorithm. We assume that the system is fully observed, i.e., the algorithm can observe $x_{t}$ at each time. We also assume that the system parameters $A$ and $B$ are known to the algorithm.

We consider the setting where a fixed horizon preview of the future cost functions and disturbances are available to the algorithm at each time step $t$. Thus, a general control policy $\pi_{t}$ at time $t$ computes the control input $u_t$ as $u_{t} = \pi_{t}(x_{1:t}, u_{1:t-1},w_{1:t+N_{w}-1}, c_{1:t+N_{c}-1})$, where $N_{w}$ and $N_{c}$ are fixed preview lengths for the disturbance and the cost functions respectively. The goal is to select a control policy $\pi = \pi_{1:T}$ in order to minimize the cumulative cost. This can be formulated as the following optimization problem
\begin{align}
    \label{eq:basic-problem-1}
    & \min_{\pi}~~  J(\pi; w_{1:T}, c_{1:T}) :=  \sum^{T}_{t = 1} c_{t}(x_{t},u_{t}), \nonumber \\
    & \text{where,}~~u_{t} = \pi_{t}(x_{1:t}, u_{1:t-1},w_{1:t+N_w-1}, c_{1:t+N_c-1}). 
\end{align}
Our goal is to characterize the performance of the optimal control policy $\pi^{*}$ and the effect of prediction horizons $N_w$ and $N_c$. More specifically, our goal is to characterize the ratio of the optimal cost to the energy in the disturbance, i.e. to characterize
\beq
    \label{eq:Dist-Gain}
 \min_{\pi}~~\sup_{w_{1:T}}~~\frac{J(\pi; w_{1:T}, c_{1:T})}{\sum^{T}_{t=1} \norm{w_{t}}^{2}}.
\eeq

We make the following assumptions on the adversarial disturbance and cost functions.
\begin{assumption}[Disturbance]
The disturbance $w_{t} \in  \mathcal{W}$ for all $t \in [T]$, where $\mathcal{W}$ is compact with $\norm{w}_2 \leq w_c$ for all $w \in \mathcal{W}$. 
\label{ass:noise}
\end{assumption}
We note that the boundedness of disturbance is a standard assumption in the online robust control literature; see \cite{agarwal2019online, simchowitz2020improper}. Let
\begin{align}
& V_t(x_t,w_{t:t+N-1}) := \underset{\tilde{u}_{t:t+N-1}}{\inf} \sum_{j = t}^{t+N-1} c_j(\tilde{x}_j, \tilde{u}_j), \nonumber \\
& ~~ \tilde{x}_{j+1} = A \tilde{x}_j + B\tilde{u}_j + w_j, ~~ \tilde{x}_t = x_t. \nonumber 
\end{align}

\begin{assumption}[System]
For all $t \in [T]$, (i) $c_{t}(\cdot, \cdot)$ is continuous and non-negative, (ii) There exists a continuous non-negative function $\sigma: \mathbb{R}^{n} \rightarrow \mathbb{R}_{+}, \underline{\alpha} > 0$ such that $c_t(x_t,u_t) \geq \underline{\alpha}\sigma(x_t)$, (iii) There exist $\overline{\alpha} > 0, \overline{\gamma} > 0$ such that $V_t(x,w_{0:N-1}) \leq \overline{\alpha}\sigma(x) + \overline{\gamma}^2 \sum_{k=0}^{N-1} \norm{w_{k}}^2_2$
\label{ass:stagecost}
\end{assumption}
We introduce a constant $\beta := \underline{\alpha}/\overline{\alpha}$. By Assumption \ref{ass:stagecost}$(ii)$ and $(iii)$, the inverse of $\beta$ defines the possible variation of the optimal cost to go for the preview $N$ at different time instances. Thus this ratio captures an important characteristic of the system, which is the variability in the cost functions. It is intuitive that along with the preview $N$, $\beta$ will also influence the gain of the optimal cost to the energy in the disturbance. Our goal is to characterize the achievable ratio of the cumulative cost to the energy in the disturbance in terms of the preview $N$ and the parameter $\beta$.

The Assumption 2.$(ii)$ states that the cost functions are lower bounded by a non-negative function. Such assumptions are typical in control settings with prediction for cost functions such as ours; see for e.g. \cite{grimm2005model}. Assumption 2.$(ii)$ is trivially satisfied by the LQR cost $c_{t}(x,u) = x^{\top}Q_{t}x + u^{\top}R_{t}u$, where $Q_{t} > 0$ and $R_{t} > 0$, with $\sigma(x)=\norm{x}^2_2$. The upper bound to $V_t$ in Assumption \ref{ass:stagecost}.$(iii)$ has two terms: the first term is the contribution of the initial state to the overall cost and the second term is the contribution of the disturbances over the period $N$ to the overall cost. Here, the factor $\overline{\gamma}^2$ is the {\it minimum achievable disturbance gain by any sequence of control actions}. The minimum achievable gain for the LQR cost function for linear systems satisfies this form of upper bound with $\sigma(x) = \norm{x}^2_2$ (see \cite[Theorem 9.5.1]{hassibi1999indefinite}). This assumption is necessary since the problem we address in this work has a solution only if a finite gain $\overline{\gamma}^2$ is achievable.



\subsection{Classical $\mathcal{H}_\infty$-Optimal Control problem}

The robust control problem we study in this work has been addressed extensively in the control community using the framework of $H_{\infty}$-optimal control. Most of these works however focus on linear systems with quadratic costs. 

Consider the system \eqref{eq:stateequation} with $c_{t}(x,u) = x^{\top}Q_{t}x + u^{\top}R_{t}u$, where $Q_{t}$ is positive semi-definite and $R_{t}$ is positive definite. Also, assume that $N_{c} = T, N_{w} = 1$. Then, the $\mathcal{H}_{\infty}$-optimal control problem can be written as
\begin{align}
    \label{eq:H-inf-formulation}
 \min_{\pi}~~\sup_{w_{1:T}}~~ \frac{J(\pi; w_{1:T}, c_{1:T})}{\sum^{T}_{t=1} \norm{w_{t}}^{2}}.
\end{align}
The goal, intuitively, is to minimize the worst-case gain from the energy in the disturbance $w_{1:T}$ to the cost incurred by the controller. In general, \eqref{eq:H-inf-formulation} is extremely hard to solve. A usual way to overcome the difficulty is to consider a relaxed problem with a given performance level $\gamma$ as given below
\begin{align}
    \label{eq:H-inf-formulation}
 \min_{\pi}~~\sup_{w_{1:T}}~~ \frac{J(\pi; w_{1:T}, c_{1:T})}{\sum^{T}_{t=1} \norm{w_{t}}^{2}} < \gamma^{2}.
\end{align}
Clearly, a solution may not exist for any arbitrary value of $\gamma$. 
If the problem is feasible for a given $\gamma$, then the optimal solution can be expressed in closed form, exploiting the quadratic nature of the cost function. The precise form of the $\mathcal{H}_{\infty}$-optimal controller can be found in \cite{doyle1988state}. 

Compared to the $\mathcal{H}_\infty$ problem, the problem we study in this work is more challenging. Our setting is an online setting in which at any point of time the controller does not know the complete future of cost functions and disturbances, but only upto a horizon $N_c < T$ and $N_w \leq N_c$. We contrast this with the classical problem, where $N_c = T, N_w = 1$; and thus it needs complete knowledge of future cost functions. Additionally, the cost functions we consider are general and need not be quadratic. Our methodology also does not need the knowledge of the optimal $\gamma$.

%% file: relatedworks.tex
\section{Related Works}
\label{sec:relatedworks}

{\it Online Stochastic Control}: This area of online control studies the online control of stochastic systems. In \cite{abbasi2011regret}, the authors study the online Linear Quadratic Regulator (LQR) problem for unknown systems and stochastic disturbances. 
The authors propose an adaptive algorithm that achieves $\sqrt{T}$ regret w.r.t the best linear feedback control policy. A drawback of their online controller is the inefficiency. In \cite{dean2018regret}, the authors propose an algorithm for the same problem that is efficient and achieves a sub-linear regret of $\mathcal{O}(T^{2/3})$. In \cite{cohen2019learning} and \cite{mania2019certainty} this result is further improved by providing an efficient algorithm with $\mathcal{O}(\sqrt{T})$. In \cite{mania2019certainty}, the authors study the same setting but with partial observations of the state and establish that $\mathcal{O}(\sqrt{T})$-regret is achievable. Recently, \cite{simchowitz2020naive} showed that $\mathcal{O}(\sqrt{T})$ is the optimal regret for the online LQR problem. 

{\it Online Robust Control}: This area of online control studies the online control of non-stochastic systems. In \cite{agarwal2019online}, the authors study the control of a known linear dynamic system with arbitrary disturbance and arbitrary controller cost function for each time step. They show that an $\mathcal{O}(\sqrt{T})$-regret is achievable with respect to the best linear feedback controller. In \cite{hazan2020nonstochastic}, the authors study the same setting but for the unknown system and show that a sub-linear regret of $\mathcal{O}(T^{2/3})$ is still achievable. Recently, \cite{simchowitz2020improper} showed that similar results are achievable with partial observation of the state for both known and unknown systems. 

{\it Receding Horizon Control}: Many receding horizon control based methods have been proposed for managing disturbances and uncertainties in the system dynamics. For example, some works handle disturbances or uncertainties by robust or chance constraints \cite{goulart2006optimization,limon2010robust, tempo2012randomized, goulart2016robust}. Adaptive RHC techniques that adapt online when the system model is unknown have also been proposed \cite{aswani2013provably,tanaskovic2019adaptive, bujarbaruah2019adaptive}. These methods primarily focus on constraint satisfaction, stability and in some cases performance improvement using the adapted models. In contrast to these works, we consider non-asymptotic performance of an online adaptive RHC. There are also considerable amount of papers that present performance analysis of RHC under both time-invariant costs \cite{angeli2011average, grune2014asymptotic, grune2015non} and time varying costs \cite{ferramosca2014economic, grune2017closed}. However, most of these studies focus on asymptotic performance. 

%% file: algorithm-results.tex
\section{Online Robust Control: Algorithm and Main Result}
\label{sec:alg-mainresults}

In this section, we present the online robust control algorithm for $N_c = N_w = N > 0$, i.e., when the preview of both the future disturbances and cost functions for period $N$ are available. Specifically, at each time $t$, algorithm has access to $c_{t:t+N-1}$ and $w_{t:t+N-1}$, in addition to the history of observation until $t$. 

The control policy we propose is a modification of the standard receding horizon controller. We denote this control policy by $\pi_{\mathrm{op}}$. Unlike the standard receding horizon control, which recomputes the decision every time step, the policy $\pi_{\mathrm{op}}$ recomputes only every $N-M+1$th time step, where $M (< N)$ will be specified later. In the period between, where the control sequences are not recomputed, the policy $\pi_{\mathrm{op}}$ selects the control input from the computed sequence for the current interval. The period $M$ overlap from one computation to the next determines how much of the preview is exploited and how frequently the control inputs are updated.

More formally, the horizon $T$ is split as the intervals $[t_1,t_2-1], [t_2,t_3-1],...$ and so on till the end of the horizon. The intervals are indexed by $i$ and the beginning of this interval is given by $t_i = (i-1)(N-M)+1$. Thus, each interval is of duration $N-M$. The control policy solves the following optimization at the beginning of interval $i$:
\begin{align}
& \underset{\tilde{u}_{t:t+N-1}}{\inf} \sum_{j = t}^{t+N-1} c_j(\tilde{x}_j, \tilde{u}_j), \nonumber \\
& \text{s.t.}~~ \tilde{x}_{j+1} = A\tilde{x}_j + B\tilde{u}_j + w_j,~~\tilde{x}_t = x_t. 
\label{eq:cont-opt-distaccess}
\end{align}
The optimization computes the optimal control sequence for the cost-to-go for the duration $N$ from the beginning of every interval. We denote the output of the optimization by $\tilde{u}_{t:t+N-1}$. Given this computed control sequence, the policy $\pi_{\mathrm{op}}$ sets the control input $u_t$ as $\tilde{u}_t$ for all $t \in [t_i, t_{i+1}-1]$. The policy $\pi_{\mathrm{op}}$ then repeats this process in every new interval. The complete algorithm for $\pi_{\mathrm{op}}$ is described in Algorithm \ref{alg:online-control-distaccess}. 

We note that policy $\pi_{\mathrm{op}}$ trivially achieves the maximum possible attenuation when $N = T$. Hence, we do not discuss this case formally. In the next theorem we characterize the disturbance gain achieved by the $\pi_{\mathrm{op}}$ when $N < T$.
\begin{algorithm}[]
\begin{algorithmic}[1]
\STATE \textbf{Input}: $N, M$
\STATE \textbf{Definition}: $t_i = (i-1)(N-M)+1$
\STATE \textbf{Initiaize}: $i = 1$
\FOR {$t = 1,...,T$}
\IF{$t = t_i$}
\STATE Observe the current state $x_t$, and the preview $c_{t:t+N-1}, w_{t:t+N-1}$ 
\STATE Solve the problem \eqref{eq:cont-opt-distaccess} to get $\tilde{u}_{t:t+N-1}$
\STATE $i = i+1$
\ENDIF
\STATE Set the control input $u_t$ as $\tilde{u}_t$ 
\ENDFOR
\end{algorithmic}
\caption{Online Robust Control with Preview ($\pi_{\mathrm{op}}$)}
\label{alg:online-control-distaccess}
\end{algorithm}

\begin{theorem}
Suppose Assumption \ref{ass:noise} and \ref{ass:stagecost} hold. Let $\pi = \pi_{\mathrm{op}}$, $M = N/2$ and $\sigma(x_1) = 0$. Let $N > 4\zeta^3$ and $\beta \geq 1/\zeta$ for some $\zeta > 1$. Then, for any $T > N$
\begin{align} 
& J(\pi_{\mathrm{op}}) \leq \gamma_{\text{op}}^2 \sum_{t=1}^T \norm{w_t}^2_2,~~ \text{where} \nonumber \\
& \gamma^2_{\text{op}} = (2\zeta+\rho(N))\overline{\gamma}^2, ~ \rho(N) \leq \mathcal{O}(1/N),  \nonumber
\end{align}
i.e., $\pi_{\mathrm{op}}$ achieves a disturbance gain $\gamma^2_{\text{op}}$.
\label{thm:orc-with-distpreview}
\end{theorem}
We discuss the analysis in Section \ref{sec:detaileddiscussion}.
\begin{remark}
The complete characterization of the bound to the cost includes an additional term that is the cost of the initial condition, $\sigma(x_1)$, as in the $\mathcal{H}_\infty$ characterization; see \cite[Theorem 9.5.1]{hassibi1999indefinite} for $\mathcal{H}_\infty$ characterization. We present the complete characterization later. In the above theorem, we assume the initial state's cost to be zero and only state the bound in terms of the disturbance gain.
\end{remark}


\begin{remark}
An implication of our general result is the following: when $\zeta = (1+o(1))$, Theorem \ref{thm:orc-with-distpreview} implies that, with preview of disturbances and cost functions $N > 4 + o(1)$, the proposed RHC achieves a disturbance gain $\gamma^2_{\text{op}} = (2+\mathcal{O}(1/N))\overline{\gamma}^2$. 
\end{remark}

%% file: analysis.tex
\section{Detail Analysis}
\label{sec:detaileddiscussion}

In the next theorem we give complete characterization of the performance of Algorithm \ref{alg:online-control-distaccess}. 
\begin{theorem}
Suppose Assumption \ref{ass:noise} and \ref{ass:stagecost} hold. Let $\pi = \pi_{\mathrm{op}}$, $N \geq 2M$ and $M > 1/\beta^2$. Then for any $T > N$,
\begin{align} 
& J(\pi_{\mathrm{op}}) \leq \mathcal{O}(\sigma(x_1)) + \omega_{\mathrm{op}}\overline{\gamma}^2 \sum_{t=1}^T \norm{w_t}^2_2,~~ \text{where} \nonumber \\
& \omega_{\mathrm{op}} = \frac{2-\beta+\kappa(M)}{\beta(1-1/(\beta^2M))}, ~ \kappa(M) = \frac{3}{\beta M} + \frac{1}{(\beta M)^2} -\frac{1}{M}. \nonumber
\end{align} 
\label{thm:orc-with-distpreview-general}
\end{theorem}
Please see the appendix 
for the detailed proof. The proof is built on the proof technique used in \cite{grimm2005model}. 

We observe that as $N$ increases, fixing $M$ to a constant does not change the $\omega_{\mathrm{op}}$. We observe this because the updates to the control actions become less frequent when $M$ is fixed and does not change with the preview $N$, which could then limit the adaptability of the control policy. It follows that we cannot set $M$ to be very small. While the optimal value for $M$ depends on the specific realization of the disturbances, $M = N/2$ is a reasonable choice considering the guarantee in Theorem \ref{thm:orc-with-distpreview-general}. In our case, when $M = N/2$, we observe that $\omega_{\mathrm{op}}$ trivially decreases as atleast $1/N$ since $\kappa$ decreases as $1/N$.

When $M = N/2$, since $\omega_{\mathrm{op}}$ is decreasing with $N$, $\omega_{\mathrm{op}} \geq \left(2-\beta\right)1/\beta$. This implies that the minimum achievable disturbance gain when $N < T$ is dependent on the system constant $\beta$. 
In addition, since $\beta \leq 1$, $\omega_{\mathrm{op}} \geq 1$, which is as expected.

The lower bound stated for $M$ is the minimum preview that is required for the response to be stable. This is because when $M$ is below this bound, the proposed policy, which in this case will only update its control actions less frequently, need not be stabilizing. Thus, we give guarantees for the problem of stabilization and disturbance gain minimization. 

\subsection{Proof of Theorem \ref{thm:orc-with-distpreview}}

\begin{proof} First, we note that $M = N/2$. Next, we note that $N > 4\zeta^3 > 2\zeta^2$. Then, given that $\zeta^2 \geq 1/\beta^2$, $N > 2/\beta^2$. That is $M > 1/\beta^2$. Thus, all the conditions in Theorem \ref{thm:orc-with-distpreview-general} are satisfied. 

Then, given that $\sigma(x_1) = 0$, $\sum_{t=1}^T  c_t(x_t,u_t) \leq \omega_{\mathrm{op}} \overline{\gamma}^2 \sum_{t=1}^T \norm{w_t}^2_2$. Substituting $M = N/2$, we get
\begin{align}
\omega_{\mathrm{op}} & = \frac{2 - \beta + \kappa(N/2)}{\beta(1 - 2/(\beta^2N))}  \nonumber \\
& \leq \frac{2 - \beta}{\beta(1 - 2/(\beta^2N))}  + \frac{ \kappa(N/2)}{\beta(1 - 2/(\beta^2N))}. \nonumber 
\end{align}

Given that $N > 4\zeta^3$ and $\beta \geq 1/(\zeta)$, we get that 
\beq
\frac{1}{2\zeta} > \frac{2 \zeta^2}{N} \geq \frac{2}{N\beta^2},~\text{i.e.,} ~ \frac{1}{\beta(1-2/(\beta^2N))} < \frac{2\zeta^2}{2\zeta -1}. \nonumber 
\eeq

By definition $\beta \leq 1$. Hence,
\beq
\kappa(N/2) = \frac{6}{\beta N} + \frac{4}{(\beta N)^2} - \frac{2}{N} = \frac{6 - 2\beta}{\beta N} +\frac{4}{(\beta N)^2} \leq \mathcal{O}(\frac{1}{N}). \nonumber 
\eeq

Thus, 
\beq 
\frac{\kappa(N/2)}{\beta(1 - 2/(\beta^2N))} \leq \mathcal{O}(\frac{1}{N}). \nonumber 
\eeq 

Then, given that $\beta \geq 1/\zeta$
\beq 
\frac{2 - \beta}{\beta(1 - 2/(\beta^2N))} \leq 2\zeta. \nonumber 
\eeq 

Thus,
\beq 
\omega_{\mathrm{op}} \leq 2\zeta + \rho(N), ~ \text{where} ~ \rho(N) \leq \mathcal{O}(1/N). \nonumber 
\eeq 

From here $\gamma^2_{\mathrm{op}}$ follows. \end{proof}

\section{Numerical Studies}
\label{sec:sim}

In this section we present numerical studies to demonstrate that the variant of the RHC algorithm we prescribe satisfies the stated bound and is better than the standard RHC approach, i.e, the RHC approach where the control inputs are recomputed every time step. We present numerical results for three examples to illustrate this. In all the examples we consider a linear dynamical system with $n = 3$ and $m = 1$. In each example, the $A$ matrix of the linear dynamical system is chosen randomly with each element of $A$ lying in the range $[0,1]$. The $B$ matrix in all the examples is given by $B = [1, 1]^\top$. The components of the disturbance lie in the range $[0,1]$. The duration of the control task in all the examples is $T = 15$ time steps. The examples differ in the cost functions and the cost functions for the respective examples are as follows,
\begin{enumerate}
    \item Example $1$ (Quadratic): the cost is quadratic, given by $c_t(x,u) = x^\top Q_tx + u^\top R_t u$, where $Q_t$ and $R_t$ are randomly chosen diagonal matrices with their diagonal elements lying in the range $[1,3]$. This example illustrates a standard time varying cost function example.
    \item Example $2$ (Non-convex): the cost is non-quadratic and is given by $c_t(x,u) = \vert (x(1)-b)\vert ^3 + (x(2)-b)^2 + u^\top u$, where $b = 0.2$. This example illustrates a non-convex cost function example. 
    \item Example $3$ (Distance to a set): the cost is given by $c_t(x,u) = a_t\sigma(x) + u^\top u$, where $a_t$ is a randomly drawn coefficient from $[0,1]$, $\sigma(x) = \norm{x}^2_\mathcal{A}$, $\norm{\cdot}_\mathcal{A}$ denotes the minimum distance to set $\mathcal{A}$. In this example, $\mathcal{A}$ is the ball of radius $0.25$ centered at $0.5$. This example is a very common objective in predictive control, where the focus is on convergence to a specific region $\mathcal{A}$ of the state space. 
\end{enumerate}
We present the numerically computed disturbance gain for two variants of the online control policy we propose: (i) $M = N-1$ (standard RHC), (ii) $M = N/2$, 

Table \ref{tab:results} shows the disturbance gain calculated by the ratio of the average total cost by the average total energy in the disturbance across $10$ iterations for both the controllers and the different cost functions. In all the examples the threshold $N$ for the guarantee in Theorem \ref{thm:orc-with-distpreview} to hold is at least $6$. Thus, we find that in all the examples the disturbance gain of the online controller with $M = N/2$ satisfies the upper bound $\gamma^2_\mathrm{op}$ given that it is less than $2\beta \overline{\gamma}^2$. We note that in the third example the theoretical bound is quite close to the actual gain suggesting the tightness of our bound. We also note that the RHC strategy we propose, given by $M = N/2$, is better than the standard RHC across the scenarios consistently.

\begin{table}\centering\resizebox{\columnwidth}{!}{
\begin{tabular}{ |c|c|c|c|c| } 
\hline
Cost & RHC Type & DG, $N = 6$ & DG, $N = 9$ & $2\beta\overline{\gamma}^2$ \\
\hline
\multirow{2}{4em}{\centering Quadratic} & Standard ($M = N-1$) & 0.585 & 0.585 & \multirow{2}{4em}{$3.68$} \\ 
& $M = N/2$ & $0.585$ & $0.585$ & \\ 
\hline
\multirow{2}{4em}{\centering Non-convex} & Standard ($M = N-1$) & $0.378$ & $0.385$ & \multirow{2}{4em}{$2.13$} \\ 
& $M = N/2$ & $0.372$ & $0.383$ & \\ 
\hline
\multirow{2}{4em}{\scriptsize Distance to a set} & Standard ($M = N-1$) & $0.384$ & $0.39$ & \multirow{2}{4em}{$0.72$} \\ 
& $M = N/2$ & $0.377$ & $0.389$ & \\
\hline
\end{tabular}}
\label{tab:results}
\caption{Disturbance Gain (DG) for different cost conditions}
\end{table}

%% file: conclusion.tex
\section{Conclusion}
\label{sec:conclusion}

In this work we study the online robust control problem with preview for linear dynamic systems with uncertainties in the cost functions and disturbance and characterize its performance in terms of the system parameters, specifically the preview length and the variability in the cost functions. Our goal is to design an online control policy that minimizes the disturbance gain, which is a well studied metric in the classical $\mathcal{H}_\infty$ framework, a standard robust control framework for quadratic costs and full preview of costs. In contrast to the $\mathcal{H}_\infty$ framework, we study the online setting where the controller only has a limited preview of the future disturbances and cost functions that is strictly less than the control horizon. The online robust controller we propose is a novel variant of the standard Receding Horizon Controller (RHC). We show that the proposed RHC can achieve a disturbance gain that differs from the minimum achievable by a function of the system dependent parameters and decreases with the preview length-$N$ as $\mathcal{O}(1/N)$. We give explicit characterization of the disturbance gain and the required threshold for the preview for the disturbance gain guarantee to hold. 
In conclusion, our work contributes towards algorithms, guarantees and mathematical techniques for the online robust control with limited preview. 


%% file: appendix.tex
\begin{appendices}



\section{Proof of Theorem \ref{thm:orc-with-distpreview-general}}
\label{sec:proof-theorem-orc-preview}
\begin{proof}
Let
\begin{align}
& \tilde{u}_{t:t+N-1} := \underset{\tilde{u}_{t:t+N-1}}{\arg\inf} \sum_{j = t}^{t+N-1} c_j(\tilde{x}_j, \tilde{u}_j), \nonumber \\
& ~~ \tilde{x}_{j+1} = A\tilde{x}_j + B\tilde{u}_j + w_j,~~\tilde{x}_t = x_t. \nonumber 
\end{align}

Let $\phi_t(k, x_t, \tilde{u}_{t:t+k-1}, w_{t:t+k-1})$ denote the state of the system at time $t+k$ following the control sequence $\tilde{u}_{t:t+k-1}$ and disturbance $w_{t:t+k-1}$. For ease of illustration, we denote $V_{t_{i}}, x_{t_{i}}, \tilde{u}_{t_i:t_i+k-1}, w_{t_i:t_i+k-1}$ and $\phi_{t_i}$ by $V_{i}, x_i, \tilde{u}^i_{0:k-1}, w^i_{0:k-1}$ and $\phi_i$. Let $n_m := N-M$, $x^j_{i+1} := \phi_{i+1}(j-1, x_{i+1}, \tilde{u}^i_{n_m:n_m+j-2}, w^{i+1}_{0:j-2})$ and $t^j_{i+1} = t_{i+1}+j-1$. Now, 
\begin{align}
& V_{i+1}(x_{i+1}, w^{i+1}_{0:N-1}) \nonumber \\
& =  \sum_{k = 0}^{N-1} c_{t_{i+1}+k}(\phi_{i+1}(k, x_{i+1}, \tilde{u}^{i+1}_{0:k-1}, w^{i+1}_{0:k-1}),\tilde{u}^{i+1}_{k})\nonumber \\
& \stackrel{(a)}{\leq} \sum_{k = 0}^{j-2} c_{t_{i+1}+k}(\phi_{i+1}(k, x_{i+1}, \tilde{u}^i_{n_m:n_m+k-1}, w^{i+1}_{0:k-1}), \tilde{u}^i_{n_m+k}) \nonumber \\
& + \inf_{\tilde{u}_{0:N-j}}\sum_{k = 0}^{N-j} c_{t^j_{i+1}+k}(\phi_{t^j_{i+1}}(k, x^j_{i+1}, \tilde{u}_{0:k-1}, w^{i+1}_{j-1:j+k-2}), \tilde{u}_k) \nonumber \\
& \stackrel{(b)}{\leq} \sum_{k = 0}^{j-2} c_{t_{i+1}+k}(\phi_{i+1}(k, x_{i+1}, \tilde{u}^i_{n_m:n_m+k-1}, w^{i+1}_{0:k-1}), \tilde{u}^i_{n_m+k})  \nonumber \\
& + \overline{\alpha} \sigma(x^j_{i+1}) + \overline{\gamma}^2 \sum_{k = t^j_{i+1}}^{t^N_{i+1}} \norm{w_k}^2. \nonumber 
\end{align}
Here, we get (a) by using the definition of $x^j_{i+1}$ and the fact that $\tilde{u}^{i+1}_{0:N-1}$ is the optimal control sequence, and (b) by applying Assumption \ref{ass:stagecost}.$(iii)$ to the second term. Similarly,
\begin{align}
& V_{i}(x_{i}, w^i_{0:N-1}) = \sum_{k = 0}^{N-1} c_{t_i+k}(\phi_i(k, x_i, \tilde{u}^i_{0:k-1}, w^i_{0:k-1}), \tilde{u}^i_{k}) \nonumber \\
& \stackrel{(c)}{=} \sum_{k=0}^{n_m-1} c_{t_i+k}(\phi_i(k, x_i, \tilde{u}^i_{0:k-1}, w^i_{0:k-1}), \tilde{u}^i_{k}) \nonumber \\
& + \sum_{k=0}^{M-1} c_{t_{i+1}+k}(\phi_{i+1}(k, x_{i+1}, \tilde{u}^i_{n_m:n_m+k-1}, \dots \nonumber \\
& w^i_{n_m:n_m+k-1}),\tilde{u}^i_{n_m+k}) \nonumber \\
& \stackrel{(d)}{=} \sum_{k=0}^{n_m-1} c_{t_i+k}(\phi_i(k, x_i, u_{t_i:t_i+k-1}, w^i_{0:k-1}), u_{t_i+k}) \nonumber \\
& + \sum_{k = 0}^{j-2} c_{t_{i+1}+k}(\phi_{i+1}(k, x_{i+1}, \tilde{u}^i_{n_m:n_m+k-1}, w^{i+1}_{0:k-1}), \tilde{u}^i_{n_m+k}) \nonumber \\
& + \sum_{k = 0}^{M-j} c_{t_{i+1}+j+k-1}(\phi_{t^j_{i+1}}(k, x^j_{i+1}, \dots \nonumber \\
& \tilde{u}^i_{n_m+j-1:n_m+j+k-2}, w^{i+1}_{j-1:j+k-2}), \tilde{u}^i_{n_m+j+k-1}). \nonumber 
\end{align}
Here, we get (c) by $\tilde{u}^i_{0:n_m-1} = u_{t_i+n_m-1}$, $t_{i+1} = t_i + n_m$ and $\phi_i(n_m, x_i, \tilde{u}^i_{0:n_m-1}, w^i_{0:n_m-1}) = x_{i+1}$, and (d) by splitting the second term further and using the definition of $x^j_{i+1}$.

Now, subtracting $V_{i}$ from $V_{i+1}$ and canceling the common terms, we get
\begin{align}
& V_{i+1}(x_{i+1}, w^{i+1}_{0:N-1}) - V_{i}(x_{i}, w^i_{0:N-1}) \leq \overline{\alpha} \sigma(x^j_{i+1}) \nonumber \\
& + \overline{\gamma}^2 \sum_{k = t^j_{i+1}}^{t^N_{i+1}} \norm{w_k}^2 \nonumber \\
& - \sum_{k=0}^{n_m-1} c_{t_i+k}(\phi_i(k, x_i, u_{t_i:t_i+k-1}, w^i_{0:k-1}), u_{t_i+k}). 
\label{eq:Vdiff-upb}
\end{align}

Let $\phi^t_i := \phi_i(t-t_i, x_i, \tilde{u}^i_{0:t-t_i-1}, w^i_{0:t-t_i-1})$. Then, from Assumption \ref{ass:stagecost}.$(ii)$, and recognizing that $x^j_{i+1} = \phi^{t_{i+1}+j-1}_i$, we get
\begin{align}
& \underline{\alpha} \sum_{t=t_i}^{t_{i+1}-1} \sigma(\phi^t_i) + \underline{\alpha} \sum_{j = 1}^{M} \sigma(x^j_{i+1}) \nonumber \\
& \leq \sum_{k = 0}^{N-1} c_{t_i+k}(\phi_i(k, x_i, \tilde{u}^i_{0:k-1}, w^i_{0:k-1}), \tilde{u}^i_{k}) \nonumber \\
& \Rightarrow ~~ \underline{\alpha} \sum_{j = 1}^{M} \sigma(x^j_{i+1}) \nonumber \\
& \leq \sum_{k = 0}^{N-1} c_{t_i+k}(\phi_i(k, x_i, \tilde{u}^i_{0:k-1}, w^i_{0:k-1}), \tilde{u}^i_{k}). \nonumber 
\end{align}

Then, by Assumption \ref{ass:stagecost}.$(iii)$, definition of $\tilde{u}^i_{0:N-1}$ and recognizing that $t^M_{i+1} = t^N_i$, we get 
\beq
\underline{\alpha} \sum_{j = 1}^{M} \sigma(x^j_{i+1}) \leq \overline{\alpha}\sigma(x_i) + \overline{\gamma}^2 \sum_{t = t_i}^{t^M_{i+1}} \norm{w_t}^2. \nonumber
\eeq

Then, there exists a $j^{*} \in [1,M]$ such that
\beq 
\sigma(x^{j^{*}}_{i+1}) \leq \frac{\overline{\alpha}}{\underline{\alpha}M}\sigma(x_i) + \frac{\overline{\gamma}^2}{\underline{\alpha}M} \sum_{t = t_i}^{t^M_{i+1}} \norm{w_t}^2. \nonumber 
\eeq 

Then, setting $j = j^{*}$ in Eq. \eqref{eq:Vdiff-upb}, we get
\begin{align}
& V_{i+1}(x_{i+1}, w^{i+1}_{0:N-1}) - V_{i}(x_{i}, w^i_{0:N-1}) \nonumber \\
& \leq \frac{\overline{\alpha}^2}{\underline{\alpha}M}\sigma(x_i)  + \frac{\overline{\alpha}\overline{\gamma}^2}{\underline{\alpha}M} \sum_{t = t_i}^{t^M_{i+1}} \norm{w_t}^2  + \overline{\gamma}^2 \sum_{t = t_{i+1}}^{t^N_{i+1}} \norm{w_t}^2 \nonumber \\
& - \sum_{k=0}^{n_m-1} c_{t_i+k}(\phi_i(k, x_i, u_{t_i:t_i+k-1}, w^i_{0:k-1}), u_{t_i+k}) \nonumber \\
& = \frac{\overline{\alpha}^2}{\underline{\alpha}M}\sigma(x_i) + \frac{\overline{\alpha}\overline{\gamma}^2}{\underline{\alpha}M} \sum_{t = t_i}^{t_{i+1}-1} \norm{w_t}^2 \nonumber \\
& + \left(\frac{\overline{\alpha}}{\underline{\alpha}M}+1\right)\overline{\gamma}^2 \sum_{t = t_{i+1}}^{t^M_{i+1}} \norm{w_t}^2 + \overline{\gamma}^2 \sum_{t = t^M_{i+1}+1}^{t^N_{i+1}} \norm{w_t}^2 \nonumber \\
& - \sum_{k=0}^{n_m-1} c_{t_i+k}(\phi_i(k, x_i, u_{t_i:t_i+k-1}, w^i_{0:k-1}), u_{t_i+k}). \nonumber 
\end{align}

Then, applying Assumption \ref{ass:stagecost}.$(ii)$ to the last term we get
\begin{align}
& V_{i+1}(x_{i+1}, w^{i+1}_{0:N-1}) - V_{i}(x_{i}, w^i_{0:N-1}) \nonumber \\
& \leq ~~\left(\frac{\overline{\alpha}^2}{\underline{\alpha}^2M}-1\right)\underline{\alpha}\sigma(x_i)  + \frac{\overline{\alpha}\overline{\gamma}^2}{\underline{\alpha}M} \sum_{t = t_i}^{t_{i+1}-1} \norm{w_t}^2 \nonumber \\
& + \left(\frac{\overline{\alpha}}{\underline{\alpha}M}+1\right)\overline{\gamma}^2 \sum_{t = t_{i+1}}^{t_{i+1}+M-1} \norm{w_t}^2 + \overline{\gamma}^2 \sum_{t = t_{i+1}+M}^{t_{i+1}+N-1} \norm{w_t}^2.  \nonumber
\end{align}

Given that $M > \frac{\overline{\alpha}^2}{\underline{\alpha}^2}$, $\frac{\overline{\alpha}^2}{\underline{\alpha}^2M} - 1 < 0$. Then, using Assumption \ref{ass:stagecost}.$(iii)$, we get
\begin{align}
& V_{i+1}(x_{i+1}, w^{i+1}_{0:N-1}) - V_{i}(x_{i}, w^i_{0:N-1}) \nonumber \\
& \leq ~~\frac{\underline{\alpha}}{\overline{\alpha}}\left(\frac{\overline{\alpha}^2}{\underline{\alpha}^2M}-1\right)\left(V_{i}(x_{i}, w^i_{0:N-1}) - \overline{\gamma}^2 \sum_{t = t_i}^{t^N_i} \norm{w_t}^2\right) \nonumber \\
& + \frac{\overline{\alpha}\overline{\gamma}^2}{\underline{\alpha}M} \sum_{t = t_i}^{t_{i+1}-1} \norm{w_t}^2 + \left(\frac{\overline{\alpha}}{\underline{\alpha}M}+1\right)\overline{\gamma}^2\sum_{t = t_{i+1}}^{t^M_{i+1}} \norm{w_t}^2 \nonumber \\
& + \overline{\gamma}^2 \sum_{t = t^M_{i+1}+1}^{t^N_{i+1}} \norm{w_t}^2 \nonumber \\
& = ~~\frac{\underline{\alpha}}{\overline{\alpha}}\left(\frac{\overline{\alpha}^2}{\underline{\alpha}^2M}-1\right)V_{i}(x_{i}, w^i_{0:N-1}) + \frac{\underline{\alpha}\overline{\gamma}^2}{\overline{\alpha}} \sum_{t = t_i}^{t_{i+1}-1} \norm{w_t}^2 \nonumber \\
& + \left(\frac{\underline{\alpha}}{\overline{\alpha}}+1\right)\overline{\gamma}^2 \sum_{t = t_{i+1}}^{t^M_{i+1}} \norm{w_t}^2 + \overline{\gamma}^2 \sum_{t = t^M_{i+1}+1}^{t^N_{i+1}} \norm{w_t}^2.  \nonumber
\end{align}

Let $a := 1 + \underline{\alpha}/\overline{\alpha}\left(\overline{\alpha}^2/(\underline{\alpha}^2M)-1\right)$. Then, it follows that $0 < a < 1$. Then
\begin{align}
& V_{i+1}(x_{i+1}, w^{i+1}_{0:N-1}) \leq aV_{i}(x_{i}, w^i_{0:N-1}) + \frac{\underline{\alpha}\overline{\gamma}^2}{\overline{\alpha}} \sum_{t = t_i}^{t_{i+1}-1} \norm{w_t}^2  \nonumber \\
& + \left(\frac{\underline{\alpha}}{\overline{\alpha}}+1\right)\overline{\gamma}^2 \sum_{t = t_{i+1}}^{t^M_{i+1}} \norm{w_t}^2 + \overline{\gamma}^2 \sum_{t = t^M_{i+1}+1}^{t^N_{i+1}} \norm{w_t}^2  \nonumber \\
& \stackrel{(e)}{\leq} aV_{i}(x_{i}, w^i_{0:N-1}) + \frac{\underline{\alpha}\overline{\gamma}^2}{\overline{\alpha}} \sum_{t = t_{i-1}+N/2}^{t_i+N/2-1} \norm{w_t}^2 \nonumber \\
& + \left(\frac{\underline{\alpha}}{\overline{\alpha}}+1\right)\overline{\gamma}^2 \sum_{t = t_{i}+N/2}^{t_{i+1}+N/2-1} \norm{w_t}^2 + \overline{\gamma}^2 \sum_{t = t_{i+1}+N/2}^{t_{i+2}+N/2-1} \norm{w_t}^2. \nonumber 
\end{align}
Here, we get (e) by lowering the initial point of the sum in second term to $t_{i-1}+N/2$, increasing the final point of the sum in last term to $t_{i+2}+N/2-1$ and then absorbing terms from the first sum and the last sum in to the middle sum by using that $M \leq N/2$. 

Then, repeatedly applying the previous inequality on the right of the inequality for each $i$, recognizing that $\sum_{t = t_{i}}^{t_{i+1}-1} c_t(x_t, u_t) \leq V_{i}(x_{i}, w^{i}_{0:N-1})$, and summing over all $i$'s, we get
\begin{align}
& \sum_{t = 1}^{T} c_t(x_t,u_t) \leq \frac{\overline{\alpha}}{1-a} \sigma(x_1) \nonumber \\
& + \frac{\left(a^2 + a(\frac{\underline{\alpha}}{\overline{\alpha}}+1)+\frac{\underline{\alpha}}{\overline{\alpha}}\right)}{1-a}\overline{\gamma}^2\sum_{t = 1}^{T} \norm{w_t}^2. \nonumber 
\end{align}

Let $\omega_{\mathrm{op}} = \frac{\left(a^2 + a(\frac{\underline{\alpha}}{\overline{\alpha}}+1)+\frac{\underline{\alpha}}{\overline{\alpha}}\right)}{1-a}$. We first observe that $a = 1+\beta(1/(\beta^2M)-1)$. Thus, $\frac{1}{1-a} = \frac{1}{\beta(1-1/(\beta^2M))}$. Substituting for $a$ in the numerator of $\omega_{\mathrm{op}}$, we get
\begin{align}
& a^2 + a(\frac{\underline{\alpha}}{\overline{\alpha}}+1)+\frac{\underline{\alpha}}{\overline{\alpha}} = a^2 + a(\beta+1)+\beta \nonumber \\
& = (1-\beta)^2 + \frac{2}{\beta M} - \frac{2}{M}+ \frac{1}{\beta^2M^2} \nonumber \\
& + (1+\beta)\left(1-\beta + \frac{1}{\beta M}\right) + \beta \nonumber \\
& = 1-2\beta + \beta^2 + \frac{2}{\beta M} - \frac{2}{M}+ \frac{1}{\beta^2M^2} + 1-\beta^2 \nonumber \\
& + \frac{1}{\beta M} + \frac{1}{M} + \beta \nonumber \\
& = 2-\beta + \frac{3}{\beta M} + \frac{1}{\beta^2M^2} - \frac{1}{M}. \nonumber 
\end{align} 

Then, letting $\kappa(M) = \frac{3}{\beta M} + \frac{1}{\beta^2M^2} - \frac{1}{M}$, we get
\beq 
\omega_{\mathrm{op}} = \frac{2 -\beta + \kappa(M)}{\beta(1-1/(\beta^2M))}. \nonumber 
\eeq 



\end{proof}

\end{appendices}